\documentclass{article}

\usepackage[utf8]{inputenc}
\usepackage[T1]{fontenc}
\usepackage{lmodern}

\usepackage{hyperref}
\usepackage{footmisc}
\usepackage{lmodern}

\usepackage{amsmath,amssymb,amsthm}

\newtheorem{theorem}{Theorem}

\newtheorem{corollary}[theorem]{Corollary}
\theoremstyle{plain}

\newcommand{\ff}[2]{\left\lfloor\frac{#1}{#2}\right\rfloor}

\newcommand{\crn}{\operatorname{cr}}


\usepackage{color}

\title{Bounding the number of non-duplicates of the $q$-side in simple drawings of $K_{p,q}$}

\author{R. Bruce Richter\thanks{rbruce@uwaterloo.ca, Department of Combinatorics \& Optimization, University of Waterloo, Ont. Waterloo, Canada} \and André C. Silva\thanks{andre.silva@ic.unicamp.br} \and Orlando Lee\thanks{lee@ic.unicamp.br, Instituto de Computação, Universidade Estadual de Campinas, Campinas - SP, 13083-852, Brazil}}

\date{}
\begin{document}

\maketitle

\begin{abstract}  The number $Z(n):=\lfloor n/2\rfloor\lfloor (n-1)/2\rfloor$ is the smallest number of crossings in a simple planar drawing of $K_{2,n}$ in which both vertices on the 2-side have the same clockwise rotation.  For two vertices $u,v$ on the $q$-side of a simple drawing of $K_{p,q}$, let $\crn_D(u,v)$ denote the total number of crossings that edges incident with $u$ have with edges incident with $v$.

We show that in any simple drawing $D$ of $K_{p,q}$ in a surface $\Sigma$ the number of pairs of vertices on the $q$-side of $K_{p,q}$ having $\crn_D(u,v)<Z(p)$ is bounded as a function of $p$ and $\Sigma$. As a consequence, we also show that, for a fixed integer $p$ and surface $\Sigma$, there exists a finite set of drawings $\mathcal{D}(p,\Sigma)$ of complete bipartite graphs such that, for each $q$, a crossing-minimal drawing of $K_{p,q}$  can be obtained by ``duplicating vertices'' in some drawing from $\mathcal D(p,\Sigma)$.

\end{abstract}

\section{Introduction}

We denote the \textit{$p$-side} and \textit{$q$-side} of the complete bipartite graph $K_{p,q}$ to be the parts of size $p$ and $q$, respectively. 

Let $u$ and $v$ be two vertices of a graph $G$ with the same neighborhood of size $p$. Let $D$ be a drawing of $G-v$ in a surface $\Sigma$ and let $\Delta$ be a neighborhood of $u$ in $\Sigma$ homeomorphic to a disk such that $\Delta$ only intersects the edges of $G$ incident with $u$. We \textit{duplicate} $u$ by drawing $v$ in the interior of $\Delta$ and drawing the edges $vw$ incident with $v$ near to $uw$ in such a way that only edges that cross $uw$ also cross $vw$. This may be done so that the edges incident with $u$ and $v$ cross at most $Z(p)=\ff{p}{2}\ff{p-1}{2}$ times in the interior of $\Delta$, as shown by Woodall~\cite{crk7n}. Note that $v$ will have the same \textit{rotation} as $u$, that is, the natural cyclic order of its neighbors induced by the counter-clockwise rotation of the edges incident with it in $D$. 


The vertex $v$ is called a \textit{duplicate} of $u$ in $D$. Any drawing $D'$ obtained from $D$ by a sequence of duplications of vertices in $D$ is called an \textit{extension} of $D$. 

A drawing of a graph $G$ in a surface $\Sigma$ is \textit{simple} if: (i) no pair of edges cross each other more than once; (ii) edges with a common incident vertex do not cross each other; and (iii) no three edges have a common crossing point. The drawing is \textit{optimal} if it has the least number of crossings over all simple drawings of $G$ in $\Sigma$. We denote the optimal crossing number of $G$ in $\Sigma$ as $\crn_\Sigma(G)$. For the sphere/plane, we omit the subscript. It is folklore that every graph has an optimal drawing that is simple. 

Duplication plays an important role in Zarankiewicz's drawings~\cite{zarankiewicz} of $K_{p,q}$ in the plane/sphere. We can obtain drawings for every $q$ starting with a planar embedding of $K_{p,2}$ and alternately duplicating the two vertices of the part of size $2$. This results in a drawing with exactly $Z(p)Z(q)$ crossings and is conjectured to be the optimal crossing number for $K_{p,q}$. This is often referred as Zarankiewicz's Conjecture. 

We denote $\crn_D(u,v)$ to be the number of crossings between the edges incident with $u$ and $v$ in a drawing $D$. Woodall~\cite{crk7n} showed that if $D$ is a simple drawing of $K_{p,q}$ on the plane/sphere and $u$ and $v$ are vertices of the $q$-side with the same rotation in $D$,  then $cr_D(u,v) \geq Z(p)$. 

With this fact we can deduce that if

\begin{equation}\label{ppty}
	cr_D(u,v) < Z(p) \text{ for every pair } u,v \text{ of the }q\text{-side} \tag{$*$}
\end{equation}
then every vertex of the $q$-side will have a unique rotation and thus $q \leq (p-1)!$. A cursory glance shows that this is not true for surfaces of higher genus, orientable or not.

In a similar fashion, our main result bounds $q$ as a function of $p$ in any drawing $D$ of $K_{p,q}$ in a surface $\Sigma$ satisfying \eqref{ppty}.

\begin{theorem}\label{thm:qbsp}
	Let $D$ be a simple drawing of $K_{p,q}$ in a surface $\Sigma$ such that, for any two vertices $u$ and $v$ of the $q$-side, $\crn_D(u,v) < Z(p)$. Then, $q$ is bounded by a function $F(p,\Sigma)$ of $\Sigma$ and $p$. 
\end{theorem}

We were only interested in the existence of the bound and do not precisely determine it. 

Note that drawings of $K_{p,q}$ that satisfy \eqref{ppty} have no pair of duplicates on the $q$-side. Drawings without duplicates are called \textit{templates}. A consequence of Theorem~\ref{thm:qbsp} is that, for a fixed integer $p$ and surface $\Sigma$, we have a finite number of templates. With that, we can derive the following result:

\begin{corollary}\label{cor:ood}
	For a fixed positive integer $p$ and surface $\Sigma$, we can obtain an optimal drawing of $K_{p,q}$ by duplicating vertices in one of finitely many drawings in a set $\mathcal{D}(p,\Sigma)$.
\end{corollary}

This is a key observation in a result of Christian, Richter and Salazar~\cite{zcffm}. They show that there exists a integer $N_0$ such that, if $\crn(K_{p,q}) = Z(p)Z(q)$ for $q \leq N_0$, then $\crn(K_{p,q}) = Z(p)Z(q)$ for all $q$. 

The proof of Theorem~\ref{thm:qbsp} follows in the next section and in Section 3 we discuss Corollary~\ref{cor:ood} and how it relates with the result of Christian~et al.

\section{Proof of Theorem \ref{thm:qbsp}}

\begin{proof}  We may assume that $p \geq 3$.

	For different vertices $i$ and $j$ of the $p$-side, if the edges $iu$ and $jv$ of $K_{p,q}$ cross each other in $D$, then there exists a $4$-cycle that self-crosses in $D$ at least once. Indeed, the edges $iu$, $jv$, $iv$, and $ju$ induce a self-crossing 4-cycle. 

	For each pair of vertices $u$ and $v$ of the $q$-side, we define a function $f_{uv}$ on the pairs of vertices $i$ and $j$ of the $p$-side such that $f_{uv}(i,j)=1$ if the $4$-cycle of $K_{p,q}$ induced by $\{i,j,u,v\}$ crosses itself in $D$, and $f_{uv}(i,j)=0$ otherwise. We note that the set of all possible such functions has size $k=2^{p \choose 2}$; therefore it is finite.

	There exists an integer $r$ such that $K_{3,r}$ is not embeddable in $\Sigma$ (see Ringel's formula for the (nonorientable) genus of complete bipartite graphs \cite{ringel65a,ringel65b}). Let $K_q$ be a complete graph such that its vertex set is the $q$-side of $K_{p,q}$. We color each edge $uv$ of $K_q$ with ``color'' $f_{uv}$. By Ramsey's Theorem, there exists a number $R:=R_k(r)$ such that if $q \geq R$, then every $k$-edge-coloring of $K_q$ with colors $1,2,\ldots,k$ contains a monochromatic copy of $K_r$. That is, the functions $f_{u,v}$ are the same for every edge $uv$ of the $K_r$; let $f$ be this common function. Note that $R$ depends on $r$ and $k$ and both depend only on $\Sigma$ and $p$, respectively. 

	Now let us define a graph $G$ whose vertex set is the $p$-side. We join $i$ and $j$ in $G$ if $f(i,j)=0$. This means that $ij \in E(G)$ if for any $u,v \in V(K_r)$ the 4-cycle induced by $\{u,v,i,j\}$ in $K_{p,q}$ does not self-cross in $D$. If there exists a triangle in $G$, then there exists a drawing of $K_{3,r}$ as a subdrawing of $D$ without crossings, which cannot happen by the choice of $r$. Thus $G$ is triangle-free. Túran's Theorem implies that $G$ has at most $\lfloor p^2/4 \rfloor$ edges. This means that there are at least ${p \choose 2} - \lfloor p^2/4 \rfloor= Z(p)$ pairs of vertices of the $p$-side that each contribute at least one crossing in $D$. Therefore, for any pair of vertices $u$ and $v$ of $K_r$, we have that $\crn_D(u,v) \geq Z(p)$. 
\end{proof}

\section{Obtaining optimal drawings}

We assert that any connected graph $G$ has only finitely many (up to homeomorphism) embeddings in a surface $\Sigma$. This fact is not new and can be implied by prior results in the literature (e.g. see Theorem 6.1 in \cite{egs}). We provide a small proof of this fact in what follows. 

An embedding is \textit{cellular} if all its faces are homeomorphic to an open disk in the plane. For a vertex $v$ of a graph $G$, a \textit{rotation} $\pi_v$ is a cyclic permutation of its incident edges. A collection of rotations for every vertex of $G$ together with a \textit{signature} function $\lambda: E(G) \rightarrow \{-1,+1\}$ composes an \textit{embedding scheme}. An embedding scheme is a combinatorial representation of an embedding of $G$. We refer the reader to \cite[Section 3.3]{gs} on how to obtain an embedding scheme from a particular embedding. For our purposes, the important thing to note here is that there are only finitely many embedding schemes for a given finite graph $G$. 

Theorem 3.1.1 in \cite{gs} proves that every cellular embedding of a graph is determined up to homeomorphism by its embedding scheme. (Furthermore, this proof includes proofs of the Jordan-Sch\"onflies Theorem and the Classification of Surfaces Theorem.)  A non-cellular embedding in a surface $\Sigma$ is obtained from a cellular embedding $\Phi’$ in a lower genus surface $\Sigma$ by adding handles and/or crosscaps to the faces of $\Phi’$.  Since each addition of a handle or crosscap increases the genus, there are only finitely many ways to convert $\Phi’$ into an embedding in $\Sigma$.  Therefore a graph $G$ has only finitely many different embeddings (cellular or not) in $\Sigma$.

We extend this result to drawings in which each pair of edges crosses at most once.  A \textit{crossing-pattern} for a graph G consists of a selection of which pairs of edges are to cross  and, for each edge, an ordering of the crossings on that edge. There are only finitely many ways this can be done (any two edges cross at most once), so G has only finitely many crossing-patterns.  Moreover, a simple drawing of G in a surface $\Sigma$ determines and is determined by an embedding of a crossing-pattern for G in $\Sigma$.  Since there are only finitely many crossing-patterns for G and only finitely many embeddings of each crossing-pattern in $\Sigma$, it follows that G has only finitely many simple drawings in $\Sigma$.


Let $\mathcal{D}(p,\Sigma)$ denote the set containing all possible simple drawings for the graphs $K_{p,q}$ with $1 \leq q \leq F(p,\Sigma)$. It is finite as per our discussion above. We now prove Corollary~\ref{cor:ood}.

\begin{proof}
%
%
%
Let $D$ be a simple drawing of $K_{p,q}$ in a surface $\Sigma$ and let $u$ and $v$ be a pair of distinct vertices of the $q$-side. If $\crn(u,v) \geq Z(p)$, then we may redraw either $u$ or $v$ as a duplicate of the other and obtain a drawing $D'$ with at most as many crossings as $D$.

Therefore some optimal drawing $D^*$ of $K_{p,q}$ exists such that, for each pair $u$ and $v$ of vertices of the $q$-side, either $cr_{D^*}(u,v) < Z(p)$ or they are duplicates. Any maximal subset $Q$ of pairwise non-duplicate vertices of the $q$-side  induces a template $B$ in $D^*$. Moreover $B$  satisfies \eqref{ppty} and, by Theorem~\ref{thm:qbsp}, is a drawing of $K_{p,q'}$ with $q' \leq F(p,\Sigma)$. In short, $D^*$ is either itself a template $B$ or an extension of some template $B$. In both cases, the template has fewer than $F(p,\Sigma)$ vertices. This means that $B \in \mathcal{D}(p,\Sigma)$. 
\end{proof}

For the plane, Christian,  Richter and Salazar~\cite{zcffm} go one step further and show that if $\crn_{B}(K_{p,q'}) = Z(p)Z(q')$, then $\crn_{D^*}(K_{p,q}) = Z(p)Z(q)$.

\hyphenation{Ci-en-tí-fi-co}
\section*{Acknowledgements}
The first author was supported by NSERC  50503-10940-500.  The second author was supported by Fundação de Amparo à Pesquisa de São Paulo Proc. 2015/04385-0, 2014/14375-9 and 2015/11937-9, Conselho Nacional de Desenvolvimento Científico e Tecnológico Proc. 311373/2015-1. The third author was financed by CNPq Proc. 303766/2018-2, CNPq Proc 425340/2016-3 and Fundação de Amparo à Pesquisa de São Paulo Proc. 2015/11937-9.

\bibliographystyle{alpha}
\bibliography{xnum}

\end{document}